\documentclass[12pt]{amsart}
\usepackage{amssymb}
\usepackage{amsmath}
\usepackage{amsthm}

\DeclareMathOperator{\im}{im}

\newcommand{\F}{{\mathcal F}}

\newcommand{\p}{\partial}

\newtheorem{lemma}{Lemma}

\newtheorem{theorem}{Theorem}

\begin{document}

\title[An invariant formula for a star product]
{An invariant formula for a star product with separation of variables}
\thanks{The paper was partially supported by the NSF grant \#1124929}
\author[Alexander Karabegov]{Alexander Karabegov}
\address[Alexander Karabegov]{Department of Mathematics, Abilene Christian University, ACU Box 28012, Abilene, TX 79699-8012}
\email{axk02d@acu.edu}

\begin{abstract} 
We give an invariant formula for a star product with separation of variables on a pseudo-K{\"a}hler manifold.
\end{abstract}
\subjclass[2010]{53D55, 53C55, 70G45}
\keywords{deformation quantization with separation of variables, K\"ahler connection, covariant tensors}

\date{July 14, 2011}
\maketitle

\section{The standard star product with separation of variables}

Given a vector space $V$, we say that the elements of the space $V[[\nu]]$ of formal series in a formal parameter $\nu$ with coefficients in $V$ are formal vectors. This way we define formal functions, formal tensors, formal operators, etc.
 
Let $M$ be a pseudo-K\"ahler manifold with the metric tensor $g_{k\bar l}$. Its inverse is a K\"ahler-Poisson tensor, $g^{\bar lk}$. The Jacobi identity for $g^{\bar lk}$ takes the form
\begin{equation}\label{E:jac}
      g^{\bar lk} \frac{\p g^{\bar qp}}{\p z^k} = g^{\bar qk} \frac{\p g^{\bar lp}}{\p z^k} \mbox{ and } 
g^{\bar lk} \frac{\p g^{\bar qp}}{\p \bar z^l} = g^{\bar lp} \frac{\p g^{\bar qk}}{\p \bar z^l},
\end{equation}
where we assume summation over repeated indices. On a contractible coordinate chart on $M$ there exists a potential $\Phi$ of the pseudo-K\"ahler metrics such that
\[
                       g_{k\bar l} = \frac{\p^2 \Phi}{\p z^k \bar z^l}.
\]
Throughout the paper we use the conventional notation
\[
             g_{k_1 \ldots k_r \bar l_1 \ldots \bar l_s} = \frac{\p^{r+s} \Phi}{\p z^{k_1} \ldots \p z^{k_r} \p \bar z^{l_1} \ldots \p \bar z^{l_s}}
\]
and the formulas
\begin{equation}\label{E:difg}
   \frac{\p g^{\bar lk}}{\p z^p} = - g^{\bar ls} g_{sp\bar t}\, g^{\bar tk} \mbox{ and } \frac{\p g^{\bar lk}}{\p \bar z^q} = - g^{\bar ls} g_{s\bar q\bar t}\, g^{\bar tk}.
\end{equation}
A star product $\ast$ on $M$ is an associative product on $C^\infty(M)[[\nu]]$ given by the formula
\begin{equation}\label{E:star}
 f\ast g = \sum_{r \geq 0} \nu^r C_r(f,g),
\end{equation}
where $C_r$ are bidifferential operators on $M$, $C_0(f,g) = fg$, and 
\[
    C_1(f,g) -C_1(g,f) = g^{\bar lk}\left( \frac{\p f}{\p \bar z^l}\frac{\p g}{\p z^k} - \frac{\p g}{\p \bar z^l}\frac{\p f}{\p z^k}\right)
\]
(see \cite{BFFLS}). We assume that the unit constant is the unity for a star product, so that $f \ast 1 = 1 \ast f = f$.
It was proved by Kontsevich in \cite{K} that deformation quantizations exist on any Poisson manifold. 

Denote by $L_f$ and $R_g$ the left star multiplication operator by $f$ and the right star multiplication operator by $g$, respectively, so that
\[
               L_f = f \ast g = R_g f.
\]
The associativity of $\ast$ is equivalent to the condition that $[L_f,R_g]=0$ for all $f,g$. A star product can be restricted to any open set in $M$.

A star product on $M$ is called a star product with separation of variables if the operators $C_r$ differentiate their first argument only in antiholomorphic directions and the second argument only in holomorphic ones. In particular,
\[
            C_1(f, g) = g^{\bar lk}\frac{\p f}{\p \bar z^l}\frac{\p g}{\p z^k}.
\]
A star product with separation of variables is characterized by the property that for a local holomorphic function $a$ and a local antiholomorphic function $b$ the operators $L_a$ and $R_b$ are the pointwise multiplication operators by the functions $a$ and $b$, respectively,
\[
      L_a = a \mbox{ and } R_b = b.
\]
A complete classification of star products with separation of variables on a pseudo-K\"ahler manifold was given in \cite{CMP1}. In particular, it was proved that on any pseudo-K\"ahler manifold $M$ there exists a global star-product with separation of variables $\ast$ uniquely determined on each contractible chart by the condition that
\[
     L_{\frac{\p \Phi}{\p z^k}} = \frac{\p \Phi}{\p z^k} + \nu \frac{\p}{\p z^k} \mbox{ and } 
R_{\frac{\p \Phi}{\p \bar z^l}} = \frac{\p \Phi}{\p \bar z^l} + \nu \frac{\p}{\p \bar z^l}.
\] 
We call it the standard star product with separation of variables. It was independently constructed in \cite{BW} using a different method. 

Explicit formulas for star products with separation of variables expressed in terms of graphs encoding bidifferential operators were given in  ~ \cite{G} and ~\cite{HX2}. In this paper we define an invariant total symbol~$T$ of the formal bidifferential operator
\begin{equation}\label{E:formbidiff}
    \sum_{r=0}^\infty \nu^r C_r
\end{equation}
\noindent which determines the standard star product with separation of variables~$\ast$ on a pseudo-K\"ahler manifold $M$. The symbol $T$ is a formal covariant tensor on $M$ separately symmetric in the holomorphic and antiholomorphic indices. The star product $\ast$ can be immediately recovered from the total symbol~$T$.  The main result of this paper is a closed invariant formula for the symbol~$T$.

\section{A closed formula for a symbol of a left star multiplication operator}

In this section we introduce a locally defined total symbol $\sigma(A)$ of a differential operator $A$ on a pseudo-K\"ahler manifold $M$. Then we give a closed formula for the total symbol $\sigma(L_f)$ of the left star multiplication operator by a function $f$ with respect to the standard deformation quantization with separation of variables on $M$.

We will work on a contractible coordinate chart on $M$ with local holomorphic coordinates $\{z^k\}$. Introduce the following locally defined operators,
\begin{equation}\label{E:def}
    \bar D^l = g^{\bar lk}\frac{\p}{\p z^k}, \ \bar D_l = \frac{\p}{\p \bar z^l} - \frac{\p^2 \Phi}{\p \bar z^l \p \bar z^q} g^{\bar qp} \frac{\p}{\p z^p} = \frac{\p}{\p \bar z^l} - \frac{\p^2 \Phi}{\p \bar z^l \p \bar z^q} \bar D^q.
\end{equation}
\begin{lemma}\label{L:canon}
  The operators
\[
    \bar z^l, \frac{\p \Phi}{\p \bar z^l}, \bar D_l, \bar D^l
\]
satisfy canonical relations. In particular, for all $l,q$,
\[
     \left [\bar D_l,\frac{\p \Phi}{\p \bar z^q}\right] = 0, [\bar D_l, \bar D^q] =0, [\bar D_l, \bar z^q] = \delta_l^q, \mbox{ and }
\left[\bar D^l, \frac{\p \Phi}{\p \bar z^q}\right] = \delta_q^l.
\]
\end{lemma}
Lemma \ref{L:canon} can be checked by direct calculations using the Jacobi identity (\ref{E:jac}) and formulas (\ref{E:difg}).

A differential operator $A$ can be uniquely represented as a finite sum of the form
\begin{equation}\label{E:normal}
       A =\sum_{s,t} \bar D_{l_1} \ldots \bar D_{l_s} f^{\bar l_1 \ldots \bar l_s}_{\bar q_1 \ldots \bar q_t}(z,\bar z) \bar D^{q_1} \ldots \bar D^{q_t}.
\end{equation}
The middle factor in each monomial in (\ref{E:normal}) is the pointwise multiplication operator by a function of $z, \bar z$. We will call this representation normal. Using auxilliary variables $\bar \eta^l$ and $\bar \zeta_l$, we introduce symbols of differential operators as functions of $z, \bar z, \bar \zeta, \bar \eta$ polynomial in $\bar \zeta$ and $\bar \eta$ such that the symbol of the operator $A$ is
\begin{equation}\label{E:symbol}
  \sigma(A) = \sum_{s,t} f^{\bar l_1 \ldots \bar l_s}_{\bar q_1 \ldots \bar q_t}(z,\bar z) \bar \zeta_{l_1} \ldots \bar \zeta_{l_s}\bar \eta^{q_1} \ldots \bar \eta^{q_t}.
\end{equation}
The composition of operators induces a composition of symbols which will be denoted by $\circ$. It is easy to check using expressions (\ref{E:normal}), (\ref{E:symbol}) and Lemma \ref{L:canon} that 
\begin{align}\label{E:prods}
\frac{\p \Phi}{\p \bar z^l} \circ F = \frac{\p \Phi}{\p \bar z^l} F,\ F \circ \frac{\p \Phi}{\p \bar z^l} = F \frac{\p \Phi}{\p \bar z^l} + \frac{\p F}{\p \bar \eta^l},\nonumber\\
    F \circ \bar\eta^l = F\bar\eta^l, \ \bar\eta^l \circ F = F\bar\eta^l  + \bar D^l F, \\
 \bar \zeta_l \circ F = F \bar \zeta_l,  \mbox{ and }  F \circ \bar \zeta_l = F \bar \zeta_l - \bar D_l F.\nonumber
\end{align}
We have the following commutators of symbols from Eqn. (\ref{E:prods}),
\begin{equation}\label{E:commut}
    \left[F,\frac{\p \Phi}{\p \bar z^l}\right]_\circ = \frac{\p F}{\p \bar \eta^l},\ [\bar\eta^l, F]_\circ  = \bar D^l F, \mbox{ and }   [\bar \zeta_l, F]_\circ = \bar D_l F.
\end{equation}
The following lemma can be obtained from (\ref{E:prods}) and (\ref{E:commut}) by standard methods.
\begin{lemma}\label{L:compos}
The composition $\circ$ of symbols $F,G$ is given by the formula
\begin{equation}\label{E:compos}
    F\circ G = \mu\left(\exp\left\{\frac{\p}{\p \bar \eta^l}\otimes \bar D^l - \bar D_l \otimes \frac{\p}{\p \bar \zeta_l} \right\}(F\otimes G)\right),
\end{equation}
where $\mu(A \otimes B) = AB$ is the pointwise product of functions. 
\end{lemma}
\noindent The exponential formal series in (\ref{E:compos}) terminates. If symbols $F,G$ do not depend on the variables $\bar \zeta$, their composition can be given by the formula
\begin{equation}\label{E:eta}
    F \circ G = \sum_{r=0}^\infty \frac{1}{r!} \frac{\p^r F}{\p \bar \eta^{l_1} \ldots \p\bar\eta^{l_r}} \bar D^{l_1}\ldots \bar D^{l_r} G.
\end{equation}

Given a formal function $f = f_0 + \nu f_1 + \ldots$, denote by $A = L_f$ the left star multiplication operator by $f$ with respect to the standard deformation quantization with separation of variables. Then $A$ is determined by the conditions that it commutes with the local operators $R_{\bar z^l}=\bar z^l$ and 
\begin{equation*}
    R_{\frac{\p \Phi}{\p \bar z^l}} =  \frac{\p \Phi}{\p \bar z^l} +  \nu\frac{\p}{\p \bar z^l}
\end{equation*}
for all $l$ and $A1 = f$. Set $F:= \sigma(A)$. Since $A$ commutes with the multiplication by the variables 
$\bar z^l$, its symbol does not depend on the variables $\bar\zeta_l$, i.e., $F = F(\nu, z, \bar z, \bar \eta)$. The condition $A1 = f$ means that 
\[
                F \big |_{\bar \eta =0} = f.
\]

The operator $R_{\frac{\p \Phi}{\p \bar z^l}}$ can be rewritten with the use of (\ref{E:def}) in the normal form as follows,
\begin{equation}\label{E:normparphi}
     R_{\frac{\p \Phi}{\p \bar z^l}} = \frac{\p \Phi}{\p \bar z^l} +  \nu \bar D_l + \nu \frac{\p^2 \Phi}{\p \bar z^l \p \bar z^q} \bar D^q.
\end{equation}
Therefore, the symbol of the operator $R_{\frac{\p \Phi}{\p \bar z^l}}$ is
\begin{equation}\label{E:symbr}
    \sigma\left(R_{\frac{\p \Phi}{\p \bar z^l}}\right) =  \frac{\p \Phi}{\p \bar z^l} +  \nu \bar \zeta_l + \nu \frac{\p^2 \Phi}{\p \bar z^l \p \bar z^q} \bar \eta^q.
\end{equation}
The condition that the symbols of $A$ and $R_{\frac{\p \Phi}{\p \bar z^l}}$ commute,
\[
    \left[F, \frac{\p \Phi}{\p \bar z^l} +  \nu \bar \zeta_l + \nu \frac{\p^2 \Phi}{\p \bar z^l \p \bar z^q} \bar \eta^q\right ]_\circ =0,
\]
can be simplified with the use of (\ref{E:commut}) and (\ref{E:eta}) to the following form:
\begin{align}\label{E:simp}
    \frac{\p F}{\p \bar \eta^l} - \nu \bar D_l F + \nu \sum_{r=0}^\infty \frac{1}{r!} \frac{\p^r F}{\p \bar \eta^{l_1} \ldots \p\bar\eta^{l_r}} \left(\bar D^{l_1}\ldots \bar D^{l_r}\frac{\p^2 \Phi}{\p \bar z^l \p \bar z^q}\right)\bar \eta^q - \\
\nonumber \nu\frac{\p^2 \Phi}{\p \bar z^l \p \bar z^q}\left(F \bar\eta^q + \bar D^q F\right) =0.  
\end{align}
Using the definition (\ref{E:def}) of the operator $\bar D_l$ we can further simplify (\ref{E:simp}):
\begin{equation}\label{E:fursimp}
     \frac{\p F}{\p \bar \eta^l} = \nu \left(\frac{\p F}{\p \bar z^l} - \sum_{r=1}^\infty \frac{1}{r!} \left(\bar D^{l_1}\ldots \bar D^{l_r}\frac{\p^2 \Phi}{\p \bar z^l \p \bar z^q}\right)\bar \eta^q \frac{\p^r F}{\p \bar \eta^{l_1} \ldots \p\bar\eta^{l_r}}\right).
\end{equation}
Introduce the following operator on symbols,
\[
     Q := \bar \eta^l\frac{\p}{\p \bar z^l} - \sum_{r=1}^\infty \frac{1}{r!} \left(\bar D^{l_1}\ldots \bar D^{l_r}\frac{\p^2 \Phi}{\p \bar z^l \p \bar z^q}\right)\bar \eta^l\bar \eta^q \frac{\p^r}{\p \bar \eta^{l_1} \ldots \p\bar\eta^{l_r}}.
\]
We see from (\ref{E:fursimp}) that the symbol $F$ of the operator $L_f$ satisfies the equation
\begin{equation}\label{E:main}
     E(F) = \nu \, Q(F),
\end{equation}
where
\[
   E = \bar\eta^l \frac{\p}{\p \bar \eta^l}
\]
is the Euler operator with respect to the variables $\bar\eta$. Denote by $\F$ the space of formal symbols which do not depend on the variables $\bar\zeta$. The elements of $\F$ are the formal symbols $S = S_0 + \nu S_1 + \nu^2 S_2 + \ldots$, where each component $S_r = S_r(z,\bar z,\bar \eta)$ is polynomial in the variables $\bar\eta$. So, in particular, $F\in \F$. The operator $E$ acts on $\F$. We have that
\[
     \F = \ker E \oplus \im E.  
\]
The kernel $\ker E$ of $E$ in $\F$ consists of the formal functions in the variables $z,\bar z$ and the image $\im E$ consists of the symbols $S = S(\nu, z, \bar z, \bar\eta)$ such that
\[
     S \big |_{\bar\eta=0} = 0. 
\]
The Euler operator $E$ is invertible on $\im E$. Denote its inverse on $\im E$ by $E^{-1}$. Observe that the operator $Q$ maps $\F$ to $\im E$. Therefore, the operator $E^{-1}Q$ is well defined on $\F$. Represent $F$ as the sum
\[
   F = f + H
\]
for some formal symbol $H$. Then $H \in \im E$ and equation (\ref{E:main}) can be rewritten as follows:
\[
    (E - \nu Q)(H) = \nu Q(f). 
\]
We see that
\[
     (1 - \nu E^{-1}Q)(H) = \nu E^{-1}Q(f).
\]
The operator $(1 - \nu E^{-1}Q)$ is invertible on $\F$. We have
\[
    F = f + H = f + \nu (1 - \nu E^{-1}Q)^{-1}E^{-1}Q(f) = (1- \nu E^{-1}Q)^{-1}(f).
\]
We have proved the following theorem:
\begin{theorem}\label{T:form}
Given a formal function $f$, the symbol $F$ of the left multiplication operator $L_f$ with respect to the standard star product with separation of variables is given by the following formula,
\[
     F = (1- \nu E^{-1}Q)^{-1}(f) = \sum_{r=0}^\infty \nu^r \left(E^{-1}Q\right)^r (f).
\]
\end{theorem}

Observe that the component $F_r$ of $F=F_0 + \nu F_1 + \nu^2 F_2 + \ldots$ is given by the formula
\[
    F_r = \left(E^{-1}Q\right)^r (f).
\]
It is the symbol of the operator $C_r(f, \cdot)$.

\section{An invariant formula for a star product with separation of variables}

Given an affine connection $\nabla$ on a manifold $M$, there exists a global vector field on the total space of the tangent bundle $\pi: TM \to M$ defined as follows. A point $y \in TM$ represents a tangent vector $v_y \in T_{\pi(y)}M$. There exists a unique vector $w_y \in T_y (TM)$ horizontal with respect to the connection $\nabla$ such that $\pi_\ast (w_y) = v_y$. Then $y \mapsto w_y$ is a global vector field on $TM$. We will denote it $w_\nabla$. In local coordinates $\{x^i\}$ on $M$ the connection $\nabla$ is determined by the Christoffel symbols $\Gamma_{ij}^k$. If $\{y^i\}$ are the fibre coordinates on $TM$ corresponding to $\{x^i\}$, then
\[
    w_\nabla = y^i \frac{\p}{\p x^i} - \Gamma_{ij}^k y^i y^j \frac{\p}{\p y^k}.
\]
A symmetric covariant tensor $S_{i_1 \ldots i_r}$ of degree $r$ on $M$ can be equivalently described as the fibrewise polynomial function
\[
    S := S_{i_1 \ldots i_r}y^{i_1} \ldots y^{i_r}
\]
of degree $r$ on $TM$. Then the function $w_\nabla S$ corresponds to the covariant tensor of degree $r+1$ which is the symmetrization of the covariant derivative of the tensor $S_{i_1 \ldots i_r}$,
\begin{equation}\label{E:sym}
     \nabla_{(i_1} S_{i_2 \ldots i_{r+1})} = \frac{1}{r+1} \sum_{\sigma \in C_{r+1}} \nabla_{i_\sigma(1)} S_{i_{\sigma(2)}\ldots i_{\sigma(r+1)}}, 
\end{equation}
where $C_{r+1}$ is the group of cyclic permutations of the set $\{1, \ldots,r+1\}$.

Introduce the following operators,
\[
    D^k = g^{lk} \frac{\p}{\p \bar z^l}.
\]
It can be shown using (\ref{E:jac}) that they pairwise commute. Consider a differential operator $A$ that commutes with the multiplication operators by the variables $\bar z$. It can be written in the normal form as a finite sum of operators of the form
\[
   A= \sum_r f_{\bar l_1\ldots \bar l_r}(z,\bar z) \bar D^{l_1}\ldots\bar D^{l_r}.
\] 
Its symbol is the following polynomial in $\bar\eta$,
\[
   \sigma(A) = \sum_r f_{\bar l_1\ldots \bar l_r}(z,\bar z) \bar \eta^{l_1}\ldots\bar \eta^{l_r}.
\]
Using Lemma \ref{L:canon}, one can recover the symbol of $A$ by the formula
\begin{equation}\label{E:con1}
   \left(\mathrm{e}^{-\bar \eta^l \frac{\p\Phi}{\p \bar z^l}} A\, \mathrm{e}^{\bar\eta^l \frac{\p\Phi}{\p \bar z^l}}\right)1 = \sum_r f_{\bar l_1\ldots \bar l_r}(z,\bar z) \bar\eta^{l_1}\ldots \bar\eta^{l_r}.  
\end{equation}

Similarly, if an operator $\tilde A$ commutes with the multiplication operators by the variables $z$, it
can be written as a finite sum
\[
   \tilde A= \sum_r \tilde f_{k_1\ldots k_r}(z,\bar z) D^{k_1}\ldots D^{k_r}.
\]
Using the fact that
\[
         \left[D^k, \frac{\p\Phi}{\p z^p}\right] = \delta_p^k
\]
one can show that
\begin{equation}\label{E:con2}
   \left(\mathrm{e}^{-\eta^k \frac{\p\Phi}{\p z^k}} \tilde A\, \mathrm{e}^{\eta^k \frac{\p\Phi}{\p z^k}}\right)1 = \sum_r \tilde f_{k_1\ldots k_r}(z,\bar z) \eta^{k_1}\ldots \eta^{k_r}.
\end{equation}
Recall that we obtain the contravariant derivative of a tensor by lifting the index of the covariant derivative,
\[
   \nabla^k = g^{\bar lk} \nabla_{\bar l} \mbox{ and } \nabla^{\bar l} = g^{\bar lk} \nabla_k.
\]
\begin{lemma}\label{L:tensors}
Given a function $f(z,\bar z)$, the following expressions,
\[
    D^{k_1}\ldots D^{k_r} f \mbox{ and } \bar D^{l_1} \ldots \bar D^{l_r} f,
\]
are symmetric contravariant tensors.
\begin{proof}
The fact that the Christoffel coefficients of the K\"ahler connection with mixed indices are all equal to zero implies that consecutively contravariantly differentiating a function $f$ in holomorphic directions we obtain the following equality,
\[
   \nabla^{k_1} \ldots \nabla^{k_r} f = D^{k_1}\ldots D^{k_r} f.
\]
This tensor is symmetric because the operators $D^k$ pairwise commute. Similarly,
\[
   \nabla^{\bar l_1} \ldots \nabla^{\bar l_r} f = \bar D^{l_1} \ldots \bar D^{l_r} f.
\]
is a symmetric contravariant tensor.
\end{proof}
\end{lemma}
The standard star product with separation of variables can be uniquely written in the form
\[
   u \ast v = \sum_{r,s} T_{k_1 \ldots k_r \bar l_1\ldots \bar l_s} \left(D^{k_1}\ldots D^{k_r}u\right) \left(\bar D^{l_1} \ldots \bar D^{l_s} v\right),
\] 
where $T_{k_1 \ldots k_r \bar l_1\ldots \bar l_s}$ is a formal covariant tensor separately symmetric in the indices $k_i$ and $\bar l_j$. Set
\[
   T := \sum_{r,s} T_{k_1 \ldots k_r \bar l_1\ldots\bar l_s} \eta^{k_1}\ldots \eta^{k_r}\bar \eta^{l_1} \ldots \bar \eta^{l_s}.
\]
We have from (\ref{E:con1}) and (\ref{E:con2}) that
\begin{equation}\label{E:conj}
  T = \mathrm{e}^{-\left(\eta^k \frac{\p\Phi}{\p z^k} +\bar \eta^l \frac{\p\Phi}{\p\bar z^l} \right)} \left(\mathrm{e}^{\eta^k \frac{\p\Phi}{\p z^k}} \ast \mathrm{e}^{\bar \eta^l \frac{\p\Phi}{\p\bar z^l}} \right).
\end{equation}
The tensor $T$ completely determines the star product $\ast$. It can be thought of as an invariant total symbol of the formal bidifferential operator (\ref{E:formbidiff}) that determines the star product $\ast$.

The operator $Q$ is a global operator on $TM$. To see it, consider the $(0,1)$-component of the vector field $w_\nabla$. Denote it by $\bar\nabla$. We have
\[
    \bar\nabla = \bar\eta^l \frac{\p}{\p \bar z^l} - \Gamma_{\bar l\bar q}^{\bar t} \bar \eta^l \bar \eta^q \frac{\p}{\p \bar \eta^t},
\]
where
\[
    \Gamma_{\bar l\bar q}^{\bar t} = g^{\bar ts} g_{s \bar l\bar q} = \bar D^{t}\frac{\p^2 \Phi}{\p \bar z^l \p \bar z^q}
\]
is the Christoffel symbol of the K\"ahler connection with antiholomorphic indices. Lowering the upper index in the curvature tensor of the K\"ahler connection, we obtain the tensor
\begin{equation}\label{E:curvlow}
    R_{kp \bar l \bar q} := g_{kp\bar n} g^{\bar nm} g_{m\bar l \bar q}  - g_{k p \bar l \bar q}.
\end{equation}
It is separately symmetric in the holomorphic and antiholomorphic indices. Lifting the holomorphic indices in (\ref{E:curvlow}) we obtain a tensor with antiholomorphic indices only,
\begin{equation}\label{E:tens}
   R_{\bar l \bar q}^{\bar l_1 \bar l_2} : = g^{\bar l_1 k_1}g^{\bar l_2 k_2} R_{k_1 k_2 \bar l \bar q}.
\end{equation}
It is separately symmetric in the lower and upper indices. It can be shown that
\[
    R_{\bar l \bar q}^{\bar l_1 \bar l_2} = - \bar D^{l_1}\bar D^{l_2}\frac{\p^2 \Phi}{\p \bar z^l \p \bar z^q}.
\]
Given $r \geq 2$, contravariantly differentiating the tensor (\ref{E:tens}) $r-2$ times in antiholomorphic directions we obtain the tensor
\begin{equation}\label{E:rlq}
  R^{\bar l_1 \ldots \bar l_r}_{\bar l \bar q} := - \bar D^{l_1}\ldots \bar D^{l_r}\frac{\p^2 \Phi}{\p \bar z^l \p \bar z^q}.
\end{equation}
The operator $Q$ can be written in an invariant form as follows:
\[
    Q = \bar\nabla + \sum_{r=2}^\infty \frac{1}{r!}R^{\bar l_1 \ldots \bar l_r}_{\bar l \bar q}\bar \eta^l\bar \eta^q \frac{\p^r}{\p \bar \eta^{l_1} \ldots \p\bar \eta^{l_r}}.
\]
It is thus globally defined on $TM$. It follows from Theorem \ref{T:form} and formula (\ref{E:con2}) that
\[
      T = \mathrm{e}^{-\eta^k \frac{\p\Phi}{\p z^k}} (1- \nu E^{-1}Q)^{-1}
\left(\mathrm{e}^{\eta^k \frac{\p \Phi}{\p z^k}} \right).
\]
Equivalently, $T$ can be obtained by applying the operator
\[
    \mathrm{e}^{-\eta^k \frac{\p\Phi}{\p z^k}} (1- \nu E^{-1}Q)^{-1}\, 
\mathrm{e}^{\eta^k \frac{\p \Phi}{\p z^k}}
\]
to the unit constant. Observing that
\[
   \mathrm{e}^{-\eta^k \frac{\p\Phi}{\p z^k}} Q \, \mathrm{e}^{\eta^k \frac{\p \Phi}{\p z^k}} = Q + \gamma,
\]
where $\gamma := g_{pq}\eta^p \bar \eta^q$, we arrive at the following theorem.
\begin{theorem}\label{T:invar}
The tensor $T$ is given by the following invariant formula:
\[
    T = \left(1 - \nu E^{-1}\left(Q + \gamma\right)\right)^{-1}1.
\]
\end{theorem}

According to \cite{HN}, around every point $x$ of a pseudo-K\"ahler manifold for any $N$ one can choose normal holomorphic coordinates such that at~$x$
\[
    g_{k_1, \ldots, k_r, \bar l} =0 \mbox{ and } g_{k,\bar l_1, \ldots, \bar l_s}=0
\]
for all $r,s \leq N$. For every $r,s \geq 2$ there exists a canonical tensor 
\[
   R_{k_1 \ldots k_r \bar l_1 \ldots \bar l_s}
\]
separately symmetric in the indices $k_i$ and $\bar l_j$ such that it coincides with $-g_{k_1 \ldots k_r \bar l_1 \ldots \bar l_s}$ at $x$ in normal coordinates around $x$ for a sufficiently large~$N$ (see \cite{HX1}). It is expressed through the tensor $g^{\bar lk}$ and covariant derivatives of the tensor (\ref{E:curvlow}). In particular,
\begin{align*}
    R_{k_1 \ldots k_r \bar l_1 \bar l_2} & = \nabla_{k_1} \ldots \nabla_{k_{r-2}} R_{k_{r-1} k_r \bar l_1 \bar l_2}\\
 \mbox{ and } & R_{k_1 k_2 \bar l_1 \ldots \bar l_s} = \nabla_{\bar l_1} \ldots \nabla_{\bar l_{s-2}} R_{k_1 k_2 \bar l_{s-1} \bar l_s}.
\end{align*}
For $r \geq 2$, the tensor $R_{k_1 \ldots k_r \bar l \bar q}$ can be obtained from the tensor (\ref{E:rlq}) by lowering the indices $\bar l_1, \ldots \bar l_r$. Similarly, for $s \geq 2$, the tensor $R_{k p \bar l_1 \ldots \bar l_s}$ can be obtained from the tensor
\[
    - D^{k_1} \ldots D^{k_s} \frac{\p^2 \Phi}{\p z^k \p z^p}
\]
by lowering the indices $k_1, \ldots k_s$. Set
\[
     \rho_{r,s} := R_{k_1 \ldots k_r \bar l_1 \ldots \bar l_s}\eta^{k_1}\ldots \eta^{k_s} \bar \eta^{l_1}\ldots \bar \eta^{l_s}.
\]
Using Theorem \ref{T:invar} and formula (\ref{E:sym}) we can easily calculate the tensor $T$ up to $\nu^4$, which allows to recover the operators $C_r$ for $r \leq 4$:
\begin{align*}
  T = 1 + \nu \gamma & + \frac{1}{2}\nu^2 \gamma^2 + \nu^3 \left(\frac{1}{6}\gamma^3 + \frac{1}{4} \rho_{2,2}\right) +\\
& \nu^4 \left(\frac{1}{24}\gamma^4 + \frac{1}{4}\gamma \rho_{2,2} + \frac{1}{12}\rho_{2,3} + \frac{1}{12}\rho_{3,2} + \frac{1}{8}\tilde \rho  \right) + \ldots,
\end{align*}
where
\[
    \tilde \rho = R_{k_1 k_2 \bar q_1 \bar q_2} g^{\bar q_1 p_1}g^{\bar q_2 p_2} R_{p_1 p_2 \bar l_1 \bar l_2}\eta^{k_1}\eta^{k_2} \bar \eta^{l_1} \bar \eta^{l_2}.
\]

\end{document}